\documentclass{article}

\usepackage{arxiv}

\usepackage[utf8]{inputenc} 
\usepackage[T1]{fontenc}    
\usepackage{hyperref}       
\usepackage{url}            
\usepackage{mathrsfs}
\usepackage{amsmath}
\usepackage{graphicx}
\usepackage{cite}
\usepackage{caption}
\usepackage{subcaption}
\usepackage{amssymb}
\usepackage{algorithm}
\usepackage[noend]{algpseudocode}

\title{On the benefit of overparameterization in state reconstruction\thanks{This research was conducted by the Australian Research Council Training Centre in Cognitive Computing for Medical Technologies (project number ICI70200030) and funded by the Australian Government.}}

\author{Jonas F.\ Haderlein,~ Iven M.\ Y.\ Mareels\thanks{I.\ Mareels is with the Centre for Applied Research, IBM A/NZ, St Leonards NSW  2065, Australia.},~ Andre Peterson,~ Parvin Zarei Eskikand,
\AND Anthony N.\ Burkitt,~ David B.\ Grayden
\thanks{J.\ Haderlein, A.\ Peterson, P. Zarei Eskikand, A.\ Burkitt and D.\ Grayden are with the Department of Biomedical Engineering, University of Melbourne VIC 3010, Australia.}
}

\begin{document}

\maketitle
\thispagestyle{empty}
\pagestyle{empty}

\begin{abstract}

The identification of states and parameters from noisy measurements of a dynamical system is of great practical significance and has received a lot of attention. Classically, this problem is expressed as optimization over a class of models. This work presents such a method, where we augment the system in such a way that there is no distinction between parameter and state reconstruction. We pose the resulting problem as a batch problem: given the model, reconstruct the state from a finite sequence of output measurements. In the case the model is linear, we derive an analytical expression for the state reconstruction given the model and the output measurements. Importantly, we estimate the state trajectory in its entirety and do not aim to estimate just an initial condition: that is, we use more degrees of freedom than strictly necessary in the optimization step. This particular approach can be reinterpreted as training of a neural network that estimates the state trajectory from available measurements. The technology associated with neural network optimization/training allows an easy extension to nonlinear models. The proposed framework is relatively easy to implement, does not depend on an informed initial guess, and provides an estimate for the state trajectory (which incorporates an estimate for the unknown parameters) over a given finite time horizon.

\end{abstract}


\section{Introduction}
Estimating states and parameters of dynamical systems given noisy output measurements is an important problem in many areas (e.g., \cite{Mareels2002, Rudy2018, Raissi2019, Moran2013, Freestone2014}) and the topic of an entire domain of research called \textit{System Identification}. We generally consider the problem of the joint estimation of states $x_t$ and parameters $\Theta$ from a known dynamical system (class of dynamical systems) modeled as
\begin{equation}
\label{statefunction}
   x_{t+1} = f(x_t, \Theta, \nu_t) , ~~  x_1 = \kappa, \\ 
\end{equation}
with measurements 
\begin{equation}
\label{measurementequation}
y_t = h(x_t, \Theta, \mu_t),
\end{equation} 
where subscript $t$ is the time index, $x_t \in R^n$, $y_t \in R^p$, $\mu_t$ represents measurement noise, $\nu_t$ is a state transition noise (say both are i.i.d and mutually independent sequences) and $\Theta \in R^m$ is an unknown parameter. The vector fields $f$ and $h$ are known and represent the model class of the state transition and measurement function, respectively. We observe that the distinction between state and parameter is somewhat artificial. Interpreting parameters as states that evolve on a different time scale, we can always rewrite the above (\ref{statefunction}) and (\ref{measurementequation}) as
\begin{equation}
\label{augmentedstatefunction}
\begin{array}{lcll}
   x_{t+1} &=& f(x_t, \theta_t, \nu_t), & x_1 = \kappa, \\ 
   \theta_{t+1} &=& \theta_t , & \theta_1 = \Theta, \\
   y_t &=& h(x_t,\theta_t,\mu_t).
   \end{array}
\end{equation}

Considering that the measurements are obtained over a limited time horizon, $y_t, t=1 \cdots N$, the problem of interest is to recover the state trajectory $x_t, t=1 \cdots N$ and the parameter $\Theta$ given the knowledge of $f, h$ and $y_t, t=1 \cdots N$. In order to address a meaningful problem we need to assume that the state and parameter are identifiable from the measurements. In broad terms, this is a generic property.

The application that motivates us stems from computational neuroscience, where we aim to estimate brain dynamics (e.g., \cite{Moran2013, Freestone2014, Hashemi2020}). Our measurements are derived from a brain computer interface, for example, the well-known Electroencephalogram (EEG). There exists substantial uncertainty about brain dynamics and system states cannot be directly measured. Brain models are usually expressed as (classes of) dynamical systems like (\ref{statefunction}). Through judicious experiments, and measurements, particular models have been constructed to capture some brain dynamics through relatively low complexity $f$ and $h$ (e.g., \cite{Meffin2004, Ipsen2020}). These models are typically used to formulate hypotheses about brain behaviour (e.g. to determine if a person is asleep, conscious, or in an epileptic condition). System identification based on $f, h$-priors using EEG-like data is of great practical interest in order to characterize meaningful system states and eventually take meaningful actions (e.g. stopping an epileptic seizure from developing). In this context, we wish to both reconstruct state trajectories and identify the parameters. Typically, the latter contain the more important information. Measurements are normally taken over a finite time window over which it may be assumed that the brain dynamics remain stationary. 

Here, we consider the problem of state and parameter estimation, from the augmented state point of view (\ref{augmentedstatefunction}) and hence assume that the dynamics are known and that we reconstruct the unknown state from given measurements and models. As a first example, we consider the special case where the known dynamics are linear. This simple case provides some intuition for the more realistic nonlinear model case. Motivated by brain dynamics, we limit the discussion to marginally stable systems (as after all, real-world measurements of brain dynamics are bounded and oscillatory in nature). 

Here, we postulate that it is advantageous to estimate the entire state trajectory $x_t, t=1 \cdots N$ and $\theta_t, t=1 \cdots N$ as one object, rather than focusing on an initial estimate $x_1, \Theta$. Our theoretical and numerical results assess the value of this `overparameterization' within the batch processing approach. 

\paragraph{State estimation approaches in control theory}
State estimation offers a set of techniques to recover a state from measurements. These methods are well studied in the literature and are typically iterative in nature, especially in the linear system case. Examples are the well-known Kalman Filter, dead-beat observer, and the Luenberger observer (e.g., \cite{Chui2017KalmanEdition, Ljung1994ModelingSystems, Anderson1990OptimalMethods}). To the best of our knowledge, state reconstruction using an `over-parameterized' loss function, as assessed in this work, has not been considered. We demonstrate feasibility and robustness against noise in the examples section and compare with these observers, acknowledging the different settings. 

\paragraph{System identification as a state estimation problem}
Reconstructing state trajectories from noisy observations can essentially be rephrased as a zero finding problem, or finding a minimum of an appropriate loss function \cite{Mareels2002}. 
This is typically a nonlinear problem, but one that is well studied and for which there are many known numerical techniques \cite{blum}. 
These methods, such as Newton's algorithm, are only guaranteed to converge when initialized sufficiently close to a minimizer and the methodology therefore depends on an initial guess or a well constructed continuation method like the homotopy ideas presented in \cite{blum}.
Although this methodology relies on the computation of the derivatives of the system in terms of its (potentially very high-dimensional) state trajectory, the approach appears tractable for many nonlinear systems \cite{Mareels2002}. Our contribution is based on these ideas and alleviates its most crucial disadvantage: the need for meaningful initial conditions. In some sense, the present approach advocates for the opposite methodology of \cite{Ortega2015}, where the
reconstruction problem is translated into the estimation of a single constant parameter which combines the initial state and model parameter.

\paragraph{Related topics in machine learning}
Several recent approaches in machine learning discuss problems that are in nature system identification problems (e.g.,\cite{Rudy2018, Raissi2019}). In \cite{Raissi2019} the authors present neural network models for dynamical systems, assuming a fraction of the unknown states are directly observed. In order to generate a direct model of the dynamics, the neural networks are trained with a dynamics-informed loss function that is similar to those presented here. In this paper, we limit our theoretical analysis to linear dynamics, but do claim that the proposed methodology is easily adapted to nonlinear cases, even with very large model classes (that is, a large dimensional $\Theta$). 

\begin{figure}[]
      \centering
      \includegraphics[scale=0.4]{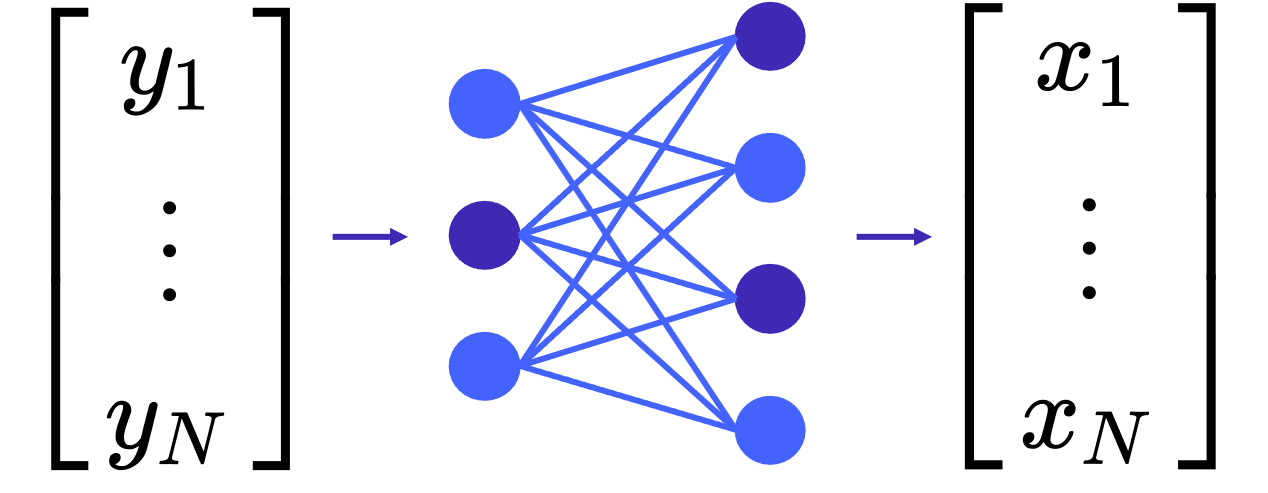}
      \caption{Graphical illustration of the estimation $\hat{X} = HY$ as represented by a dense neural network. Here, $H$ can be thought of as the weights of the network and $Y$ as the input.}
      \label{figure1}
\end{figure}

The paper is organized as follows. To set the scene, we start with known linear dynamics and present a batch filter to simultaneously estimate the entire state trajectory $x_t, t= 1 \cdots N$ via the minimization of a loss function as in \cite{Mareels2002}. In the linear case, the resulting estimator has an analytical expression. This estimator is then viewed as a linear neural network that takes noisy measurements as input, and outputs the state trajectory (see Fig. \ref{figure1}). 
From this, we then step into a nonlinear state and parameter reconstruction problem, using such a neural network approach and similar loss function. 
The ensuing nonlinear optimization problem is solved via gradient descent methods, using the automatic
differentiation methods inherent in typical deep learning toolboxes (e.g., \cite{Abadi2016TensorFlow:Learning}).

It transpires that generically a scalar valued series of observations (of sufficient length as dictated by the state dimension) suffices. The genericity of observability and the minimum window length follow from the Takens-Aeyels-Sauer theorem \cite{Takens1981}\footnote{This theorem states that state reconstruction from output measurements is inter alia a generic property of a large class of dynamical systems.}.

\section{Preliminaries}
As flagged, we start with the state (trajectory) estimation for a linear, known, finite dimensional and time invariant system when both the state transition and the measurement are affected by additive white noise\footnote{We only present this case here, but the below methodology can easily be adapted to time varying systems with known matrices $A_t$ as well as non-autonomous systems $x_{t+1}=A x_t + Bu_t + \nu_t$ with known $B$ and $u_t$.}. In particular,
\begin{equation}
\label{linearwithnoise}
x_{t+1}=A x_t + \nu_t,~~ y_t = C x_t + \mu_t,
\end{equation}
where $\nu_t$ is an i.i.d. sequence of normal random variables with zero mean and positive definite covariance matrix $\alpha >0$ and $\mu_t$ is equally an i.i.d. sequence of normal random variables with zero mean and variance $\beta>0$. We assume $\nu$ and $\mu$ to be independent random variables. In the first instance, we have knowledge of both the transition matrix $A \in R^{n \times n}$ as well as the output matrix $C \in R^{1 \times n}$. We also assume that the pair $(A,C)$ is observable\footnote{This is a generic property but, in a state estimation problem, it is without loss of generality as, without observability, the problem is not well posed.}.

Furthermore, we assume that we have measured a sequence of outputs $Y = (y_1 \cdots y_N)^T$, where $N>n-1$. 
Given the prior information ($A, C, Y$), we estimate the state trajectory $x_t$ over the finite horizon $t=1 \cdots N$, or the vector $X=({x}_1^T \cdots {x}_N^T)^T$.  

To find an appropriate estimate $\hat{X}$, we propose to minimize a non-negative loss function $\ell$ that incorporates all our prior knowledge:
\begin{equation}
\label{minimization}
     \min_{\hat{X}} \ell(\hat{X}|A,C,Y).
\end{equation} 
We discuss the loss function in Section \ref{mainsection}. For this loss, we show that an estimate $\hat{X}$ can be derived analytically in the linear case. In case the model is not linear, an algorithmic approach is required that maps all available information to an estimate $\hat{X}$. In general, the state estimate $\hat{X}$ is a function of the available output measurements and the chosen loss function that incorporates the prior information:
\begin{equation} 
    \hat{X} = F( Y, \ell).
\end{equation} 
For example, $F$ may take the form of a recipe involving the optimization of the chosen loss function through an artificial neural network approach. 
For the linear case, we may postulate that our estimate for $x_t$ is the output of a linear neural network of the form $\hat{x}_t =H_td$ for some appropriate $H_t \in R^{n \times \delta}$, for all $t=1 \cdots N$, and a particular non-zero input vector $d \in R^{\delta}$, or
\begin{equation}
\label{network}
    \hat{X}=Hd,
\end{equation}
with $H =(H_1^T \cdots H_N^T)^T$. 

A special case of the above is to select $d=Y$, or $\hat{X}=HY$, which may be of special interest since we find a direct 
mapping from $Y$ to $X$ quite appealing\footnote{In case of a neural network, e.g. Neural Tangent Kernel (NTK) parameterization corrects for different sizes of $d$ and thus $H$ \cite{Jacot2018}.}.


\section{State Reconstruction via Overparameterized Optimization}
\label{mainsection}

First, consider the linear system without noise: $x_{t+1}=Ax_t$ and $y_t=Cx_t$, and $t= 1 \cdots N$. To find an appropriate $\hat{X}$, minimize the following non-negative loss function $\ell$:
\begin{equation}
\label{loss}
    \ell(\hat{X}) = \sum_{t=1}^{N-1} \| \hat{x}_{t+1}-A\hat{x}_t\|^2+ \rho \sum_{t=1}^{N}\|y_t-C\hat{x}_t\|^2 \geq 0.
\end{equation}
Here, $\rho>0$ reflects the reliability of the model versus the measurements; the larger $\rho$, the more trust we place on the measurements. Observe that in this noise free case $\hat{x}_t=x_t$, for $t=1 \cdots N$, the loss would indeed be zero. 

Overparameterization is manifestly present in that the estimate of the complete state vector $\hat{X}$ is pursued and not just the estimation of the initial condition $\hat{x}_1$, from which a state estimate could follow via $\hat{x}_t=A^{t-1}\hat{x}_1$.

As may be expected in a linear-quadratic setting, there is a unique, optimal $\hat{X}$. To show this, let us introduce the following matrices:
\begin{equation}
    \mathscr{A}= \left( \begin{array}{rrrcrr}
    -A & I & 0 & \cdots & 0 & 0 \\
    0 & -A & I & \cdots & 0 & 0 \\
     & & & \cdots & & \\
     0 & 0 & 0 & \cdots & -A & I \end{array} \right) \in R^{n(N-1) \times nN}
\end{equation}
and 
\begin{equation}
    \mathscr{C}= \left( \begin{array}{cccc}
    C & 0_n & \cdots & 0_n \\
    0_n & C & \cdots & 0_n  \\
     & &  \cdots &  \\
     0_n & 0_n & \cdots & C \end{array} \right) \in R^{N \times nN},
\end{equation}
where $I$ is the identity matrix of the same dimension as $A$, $0$ is a matrix of the same dimension as $A$ with all elements zero, and $0_n$ is a row vector of the same dimension as $C$ with zero elements. 

Using the matrices $\mathscr{A}$ and $\mathscr{C}$, the loss function (\ref{loss}) may be rewritten as
\begin{equation}
    \ell(\hat{X})=\|\mathscr{A}\hat{X}\|^2 + \rho \|Y-\mathscr{C}\hat{X}\|^2.
\end{equation}

\subsection{Uniqueness}
Using the inner product form of the Euclidean norm, the loss function can be written
\begin{equation}
\ell(\hat{X})= \hat{X}^T\mathscr{A}^T\mathscr{A}\hat{X} + \rho (Y-\mathscr{C}\hat{X})^T(Y-\mathscr{C}\hat{X}).
\end{equation}
This is a scalar valued, quadratic expression in the unknown $\hat{X}$. 

The gradient of $\ell$ with respect to $\hat{X}$ is given by
\begin{equation}
    \nabla \ell(\hat{X}) = 2\left( \mathscr{A}^T\mathscr{A} + \rho \mathscr{C}^T\mathscr{C} \right) \hat{X} - 2 \rho \mathscr{C}^TY.
\end{equation}
Setting this gradient equal to zero, we obtain a linear equation that we can solve for $\hat{X}$,
\begin{equation}
    \left( \mathscr{A}^T\mathscr{A} + \rho \mathscr{C}^T\mathscr{C} \right) \hat{X} = \rho \mathscr{C}^TY.
\end{equation}

Because the pair $(A,C)$ is observable, there is a unique solution $\hat{X}^*$, which, with $\mathscr{O} = \mathscr{A}^T\mathscr{A}+ \rho \mathscr{C}^T\mathscr{C}$, is given by
\begin{equation}
\label{linearestimator}
    \hat{X}^* =  \mathscr{O}^{-1} \rho \mathscr{C}^TY.
\end{equation}
That $\mathscr{O}$ is nonsingular if and only if $(A,C)$ is observable may be argued as follows. If the matrix $\mathscr{O}$ were singular, then there exists a vector $z \neq 0$ such that $(\mathscr{A}^T\mathscr{A} + \rho \mathscr{C}^T\mathscr{C}) z = 0$ and, because of symmetry, both $\mathscr{A}z=0$ and $\mathscr{C}z=0$. Observability implies that this is only possible for $z=0$, which contradicts the starting premise. 

Given uniqueness, and that $X$ itself minimizes the loss function by construction, it also follows that $\hat{X}=X$ and that $\ell(\hat{X}^*)=0$.

\subsection{Noise Analysis}

Consider now the estimation when both the state transition and the measurement are affected by additive white noise. In particular, we consider the dynamics to be described as in (\ref{linearwithnoise}).

Introduce the notation $\nu = ( \nu_1^T \cdots \nu_N^T)^T$ and $\mu= (\mu_1 \cdots \mu_N)^T $.


The estimate follows as per the above (\ref{linearestimator}). Given $Y$ is affected by both the measurement noise $\mu$ as well as the state transition noise $\nu$, the estimate $\hat{X}$ is a random variable, which will differ from the actual state $X$ that only depends on the state transition noise $\nu$. By optimality of the estimate, it must be that $\ell(\hat{X}^*) \le \ell(X)$. The difference can be computed as follows.

The cost for the actual state trajectory $X$ is given by
\begin{equation}
\ell(X) = (X-\mathscr{O}^{-1}\rho \mathscr{C}^TY)^T\mathscr{O}(X-\mathscr{O}^{-1}\rho \mathscr{C}^TY) 
+ \rho Y^TY - \rho^2 Y^T\mathscr{C}\mathscr{O}^{-1}\mathscr{C}^TY.
\end{equation}
The minimum cost is attained by the estimate (\ref{linearestimator}) and is given by
\begin{equation}
    \ell(\hat{X}^*) = \rho Y^TY - \rho^2 Y^T\mathscr{C}\mathscr{O}^{-1}\mathscr{C}^TY .
\end{equation}

To analyze the difference between the estimate $\hat{X}$ and the actual state $X$ in more detail, and to consider the effect of the noises, we introduce the initial condition response sequence,
\begin{equation}
x_{i,t+1}=A x_{i,t},~~ y_{i,t} = C x_{i,t}, ~~x_{i,1}=x_1,
\end{equation}
and a zero-initial condition, noise response sequence,
\begin{equation}
x_{n,t+1}=A x_{n,t}+\nu_t,~~y_{n,t} = C x_{n,t}+\mu_t, ~~ x_{n,1}=0.
\end{equation}

Given linearity of the dynamics, we have $x_t=x_{i,t}+x_{n,t}$ and $y_t=y_{i,t}+y_{n,t}$. Let $X_n=(x_{n,1}^T \cdots x_{n,N}^T)^T$, $X_i=(x_{i,1}^T \cdots x_{i,N}^T)^T$, and similarly $Y_n=(y_{n,1} \cdots y_{n,N})^T$ and $Y_i=(y_{i,1} \cdots y_{i,N})^T$. Hence, $X=X_n+X_i$ and $Y=Y_n+Y_i$.

We now reconsider the loss function and rewrite it as
\begin{equation}
\ell(\hat{X})=
(\hat{X}-\mathscr{O}^{-1}\rho \mathscr{C}^T(Y_n+Y_i))^T\mathscr{O}(\hat{X}-\mathscr{O}^{-1}\rho \mathscr{C}^T(Y_n+Y_i))
+ (Y_n+Y_i)^T (\rho I-\rho^2 \mathscr{C}\mathscr{O}^{-1}\mathscr{C}^T)(Y_n+Y_i).
\end{equation}
Further, with $X_i=\mathscr{O}^{-1}\rho \mathscr{C}^TY_i$,
$Y_n=\mathscr{C}X_n+\mu$, and recognizing that $(I-\rho \mathscr{C}\mathscr{O}^{-1}\mathscr{C}^T)Y_i=0$, we find the loss value for $X$ to be
\begin{equation}
   \ell(X)=
   (X_n-\mathscr{O}^{-1}\rho \mathscr{C}^T(\mathscr{C}X_n+\mu))^T\mathscr{O}(X_n-\mathscr{O}^{-1}\rho
   \mathscr{C}^T(\mathscr{C}X_n+\mu))
   +Y_n^T(\rho I-\rho^2 \mathscr{C}\mathscr{O}^{-1}\mathscr{C}^T)Y_n.
\end{equation}

For the optimal estimate $\hat{X}^*=\mathscr{O}^{-1}\rho \mathscr{C}^TY$,
we have the loss
\begin{equation}
       \ell(\hat{X}^*)= Y_n^T(\rho I-\rho^2 \mathscr{C}\mathscr{O}^{-1}\mathscr{C}^T)Y_n.
\end{equation}
Clearly, we can conclude that $\ell(\hat{X}^*) \le \ell(X)$. 

To probe further, we may also consider the difference in expected cost over all noise sequences, which leads to
\begin{multline}
   E\ell(X)=
   EX_n^T(I-\rho \mathscr{O}^{-1}\mathscr{C}^T\mathscr{C})^T\mathscr{O}(I-\rho \mathscr{O}^{-1}\mathscr{C}^T\mathscr{C})X_n \\
   + E \rho^2 \mu^T\mathscr{C}\mathscr{O}^{-1}\mathscr{C}^T\mu 
   + EX_n^T\mathscr{C}^T(\rho I-\rho^2 \mathscr{C}\mathscr{O}^{-1}\mathscr{C}^T)\mathscr{C}X_n 
   + E\mu^T(\rho I-\rho^2 \mathscr{C}\mathscr{O}^{-1}\mathscr{C}^T)\mu,
\end{multline}
or
\begin{equation}
   E\ell(X)=EX_n^T\mathscr{A}^T\mathscr{A}X_n + \rho E\mu^T\mu.
\end{equation}
Whereas, for the optimal cost, we have 
\begin{equation}
  E\ell(\hat{X}^*)= EX_n^T\mathscr{C}^T(\rho I-\rho^2 \mathscr{C}\mathscr{O}^{-1}\mathscr{C}^T)\mathscr{C}X_n
   +   E\mu^T(\rho I-\rho^2 \mathscr{C}\mathscr{O}^{-1}\mathscr{C}^T)\mu,
\end{equation}
or
\begin{equation}
   E\ell(\hat{X}^*)= EX_n^T\mathscr{A}^T\mathscr{A}X_n-EX_n^T\mathscr{A}^T\mathscr{A}\mathscr{O}^{-1}\mathscr{A}^T\mathscr{A}X_n
   + E\mu^T(\rho I-\rho^2 \mathscr{C}\mathscr{O}^{-1}\mathscr{C}^T)\mu.
\end{equation}
Hence, it is clear that
\begin{equation}
\label{expectednoise}
    E\ell(X) = E\ell(\hat{X}^*)+EX_n^T\mathscr{A}^T\mathscr{A}\mathscr{O}^{-1}\mathscr{A}^T\mathscr{A}X_n
    +\rho^2 E\mu^T\mathscr{C}\mathscr{O}^{-1}\mathscr{C}^T\mu,
\end{equation}
and
\begin{equation}
    E\ell(X) \ge E\ell(\hat{X}^*),
\end{equation}
with equality if and only if the noise terms are zero. 

The expression (\ref{expectednoise}) also shows how our choice of $\rho$ affects the results.


\section{Examples}
\label{results}

\subsection{Example 1: $A=C=1$}

\begin{figure}[]
     \centering
     \begin{subfigure}[]{3.1cm}
         \centering
         \includegraphics[height = 3.1cm]{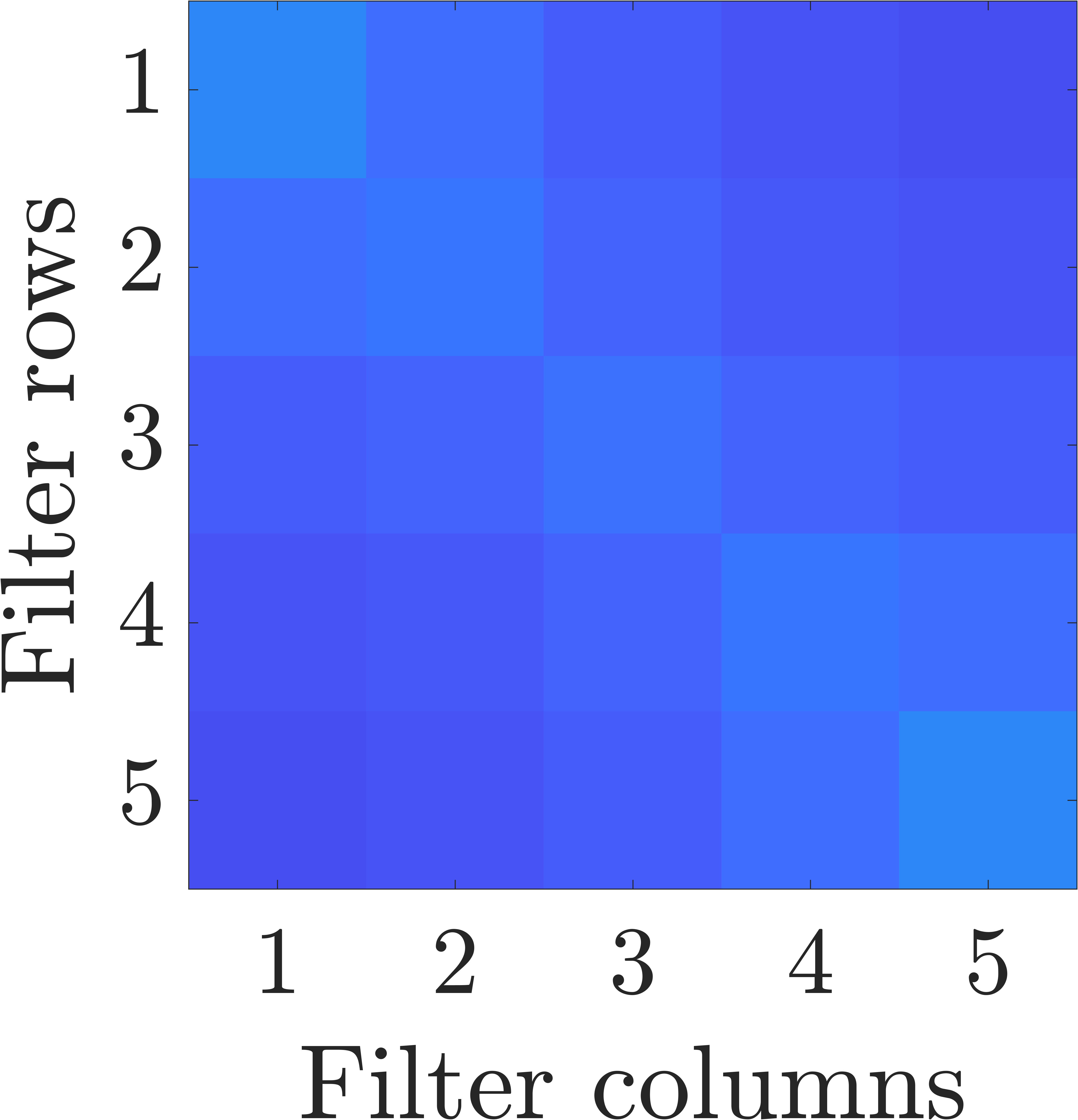}
         \caption{}
         \label{}
     \end{subfigure}
     \begin{subfigure}[]{3.1cm}
         \centering
         \includegraphics[height = 3.1cm]{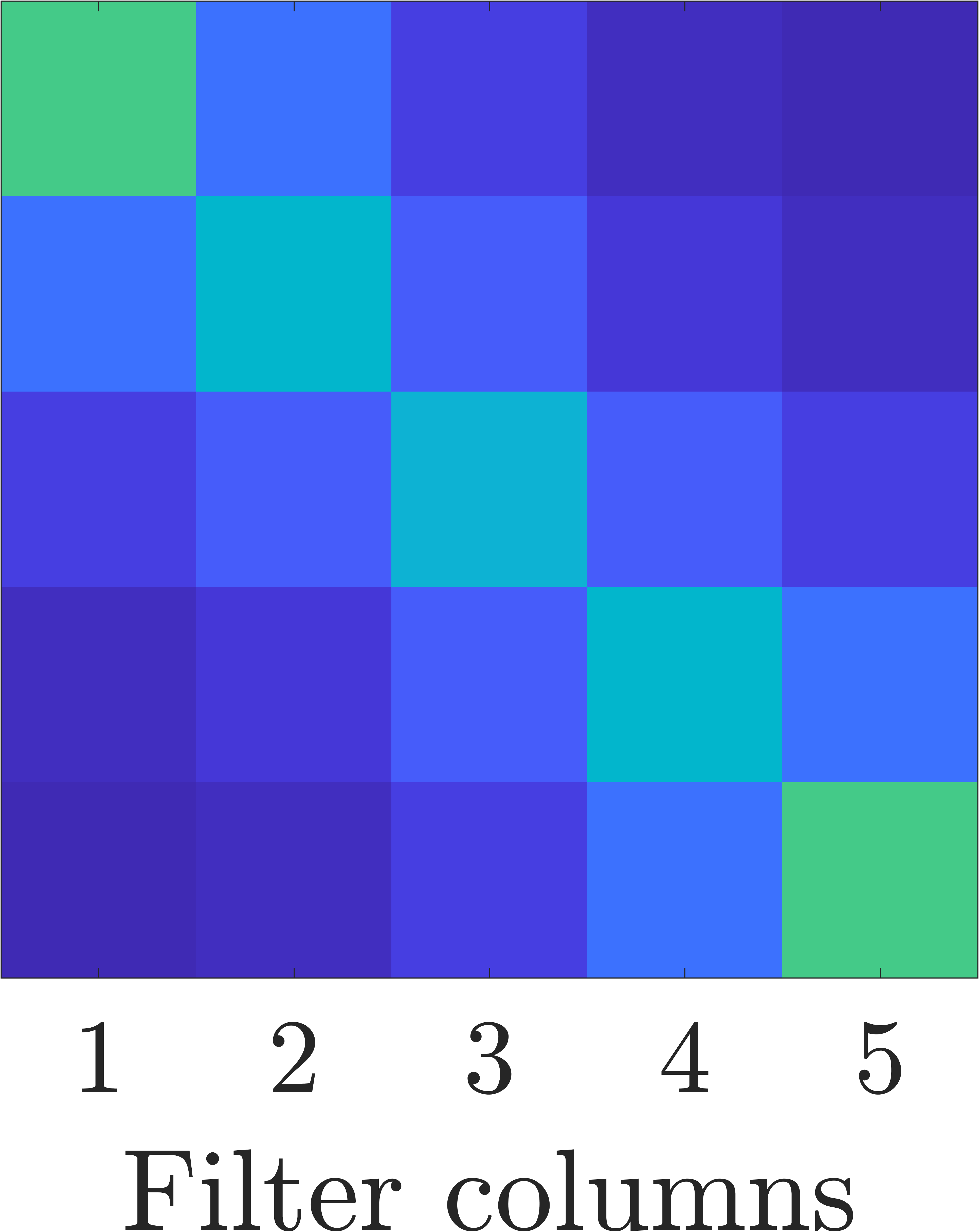}
         \caption{}
         \label{}
     \end{subfigure}
     \begin{subfigure}[]{3.1cm}
         \includegraphics[height = 3.1cm]{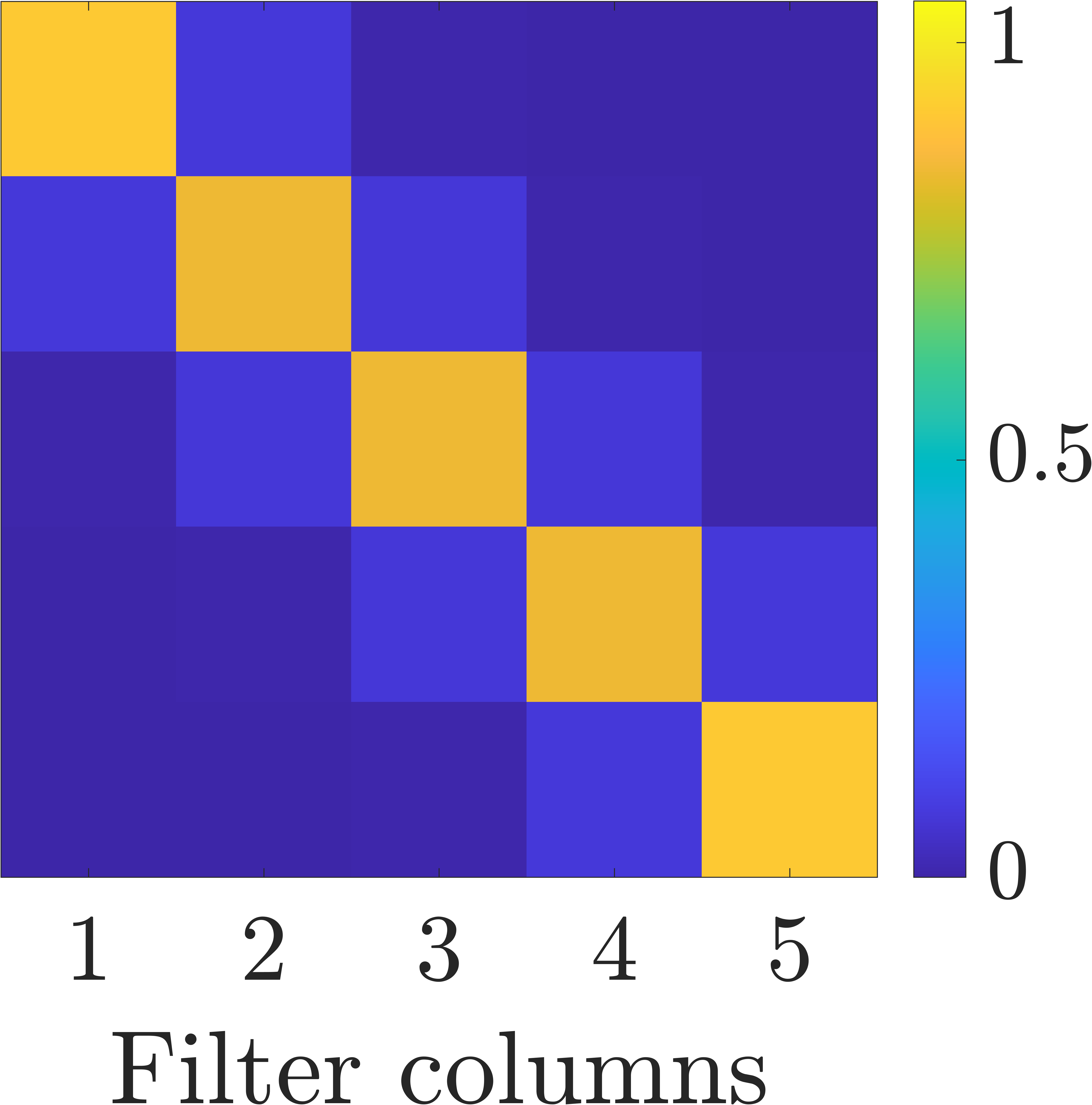}
         \caption{}
         \label{}
     \end{subfigure}
        \caption{Filter matrices $H^*=\mathscr{O}^{-1} \rho \mathscr{C}^T$ for  $A=C=1$ and different $\rho$: a) $\rho=0.1$, b) $\rho=1$, c) $\rho=10$. Colors correspond to weights as indicated by the colorbar.}
        \label{fig:Filters}
\end{figure}

To gain some intuition about the role of $H^* = \mathscr{O}^{-1} \rho \mathscr{C}^T$, let us consider the simple scalar dynamics case $A=C=1$ with $N=5$. Fig. \ref{fig:Filters} illustrates $\hat{X}=H^*Y$ for different weightings $\rho>0$.

Observe the importance of using all measurements for every time instant state estimate for small $\rho$ (Fig.\ \ref{fig:Filters}(a)), as well as the total reliance on a single observation $y_t$ as $\rho$ grows (Fig.\ \ref{fig:Filters}(c)).

The rows in $H^* = \mathscr{O}^{-1} \rho \mathscr{C}^T$ can be interpreted as filters that take $Y$ to a single $\hat{x}_.$. As $\rho$ shrinks and $N$ grows, this filter approximates more and more the least squares estimate of a constant given noisy measurements.  

\subsection{Example 2: Comparison with dead-beat observers}
\label{example2}

Consider a system $x_{t+1}=A x_t + \nu_t,~ y_t = C x_t + \mu_t$ with $n=10$, where $A$ is a companion matrix with 5 complex conjugate pairs of eigenvalues of magnitude 1 (the angles of the eigenvalues above the real axis are $\pi/2 + \pi/10 i,~i=0,1,\cdots,4 $), $C=[1,0, \cdots,0]$, $x_1=[1,1,\cdots,1]^T$, and both state and measurement noise are i.i.d. with a diagonal covariance matrix $\alpha$ with $\sigma_{\nu}^2$ on the diagonal, and measurement noise variance $\beta = \sigma_{\mu}^2$.

We compare the estimate $\hat{X} = \mathscr{O}^{-1} \rho \mathscr{C}^TY$ with the state estimate of a simple dead-beat observer and a dead-beat observer with sliding window for different values of $N$ and $\sigma_{\nu}$ (see Algorithm \ref{alg:the_alg}). To this end, we simulate $100$ state and measurement sequences for each combination of $N$ and $\sigma_{\nu}$, and compute the mean of the relative error $e = \|\hat{X}-X\|/\|X\|$ for each method. We selected $\sigma_{\mu} = 10 \sigma_{\nu}$, and set $\rho = 1$.

\def\NoNumber#1{{\def\alglinenumber##1{}\State #1}\addtocounter{ALG@line}{-1}}

\begin{algorithm}[H]
	\caption{Noise analysis (given $A, C, \rho$)}
	\label{alg:the_alg}
	\begin{algorithmic}[1]
	\For {$\sigma_{\nu} \in [0,0.1,\ldots, 1]$ and $N \in [10, 20, \ldots, 100]$}
		\State Simulate ground-truth with $x_1=[1,1,\ldots,1]^T$,
		\NoNumber{$x_{t+1}=A x_t + \nu_t,~ y_t = C x_t + \mu_t$, $t=1,\ldots,N,$}
		\NoNumber{$\nu_t \sim N(0,\sigma_{\nu}^2)$, $\mu_t \sim N(0,100\sigma_{\nu}^2)$, $\forall t$}
		\State Concatenate $X, Y$
		
        \State Calculate $\mathscr{A}$, $\mathscr{C}$, and $\mathscr{O}^{-1}$
        \State $\hat{X}_H = \mathscr{O}^{-1} \rho \mathscr{C}^TY$
        \State $e_H = \|\hat{X}_H-X\|/\|X\|$ \hfill \textit{(proposed method)}
        
        \State $O_{DB} = [C^T, (CA)^T, \ldots (CA^{N-1})^T]^T$
        \State $\hat{X}_{DB} = [I^T, A^T , \ldots, (A^{N-1})^T]^T O_{DB}^+ Y$
        \State $e_{DB} = \|\hat{X}_{DB}-X\|/\|X\|$ \hfill \textit{(simple dead-beat obs.)}
        
        \State $O_{sDB} = [C^T, (CA)^T, \ldots, (CA^{n-1})^T]^T$
        \For {$t=1:N-n+1$}
            \State $\hat{x}_{t,sDB} = O_{sDB}^{+}[y_t,\ldots,y_{t+n-1}]^T$
        \EndFor
        \For {$t=N-n+2:N$}
            \State $\hat{x}_{t,sDB} = A\hat{x}_{t-1}$
        \EndFor
        \State Concatenate $\hat{X}_{sDB}$
        \State $e_{sDB}=\|\hat{X}_{sDB}-X\|/\|X\|$\hfill\textit{(sliding dead-beat obs.)}
        
		\EndFor
	\end{algorithmic} 
\end{algorithm}

\begin{figure}
     \centering
     \begin{subfigure}[b]{0.55\textwidth}
         \centering
         \includegraphics[width=\textwidth]{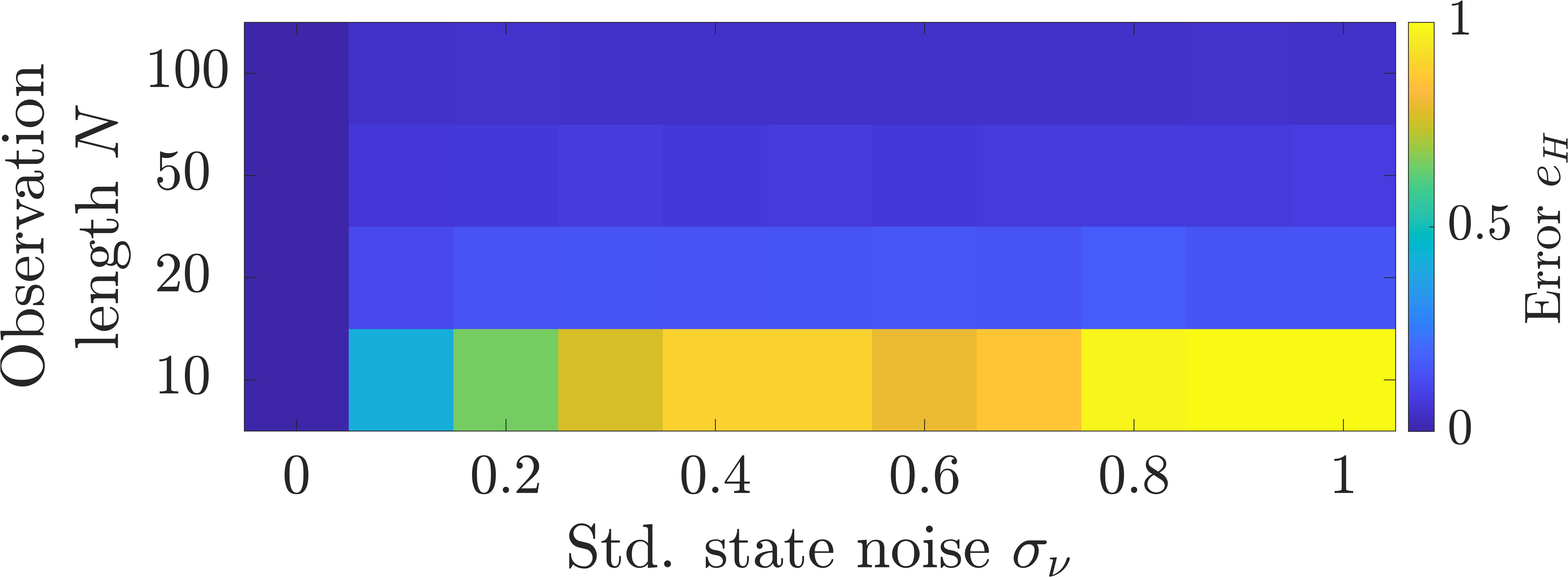}
         \caption{}
         \label{}
     \end{subfigure}
     \begin{subfigure}[b]{0.55\textwidth}
         \centering
         \includegraphics[width=\textwidth]{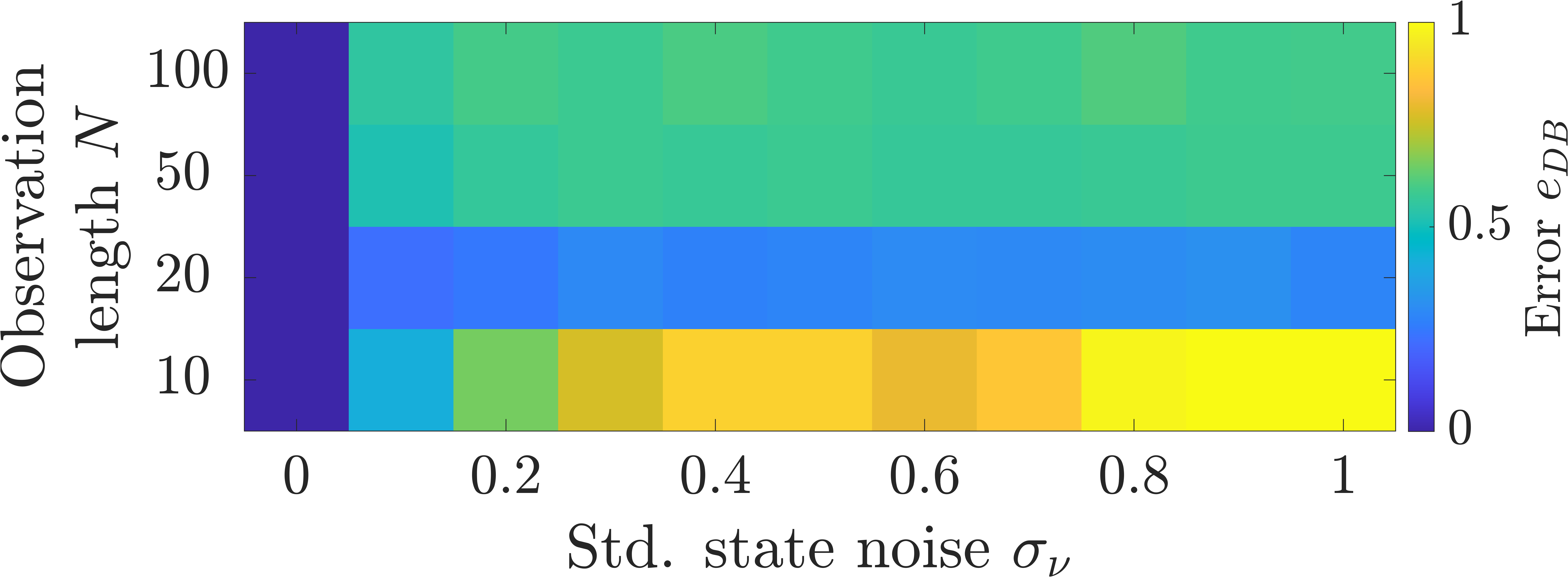}
         \caption{}
         \label{}
     \end{subfigure}
     \hfill
     \begin{subfigure}[b]{0.55\textwidth}
         \centering
         \includegraphics[width=\textwidth]{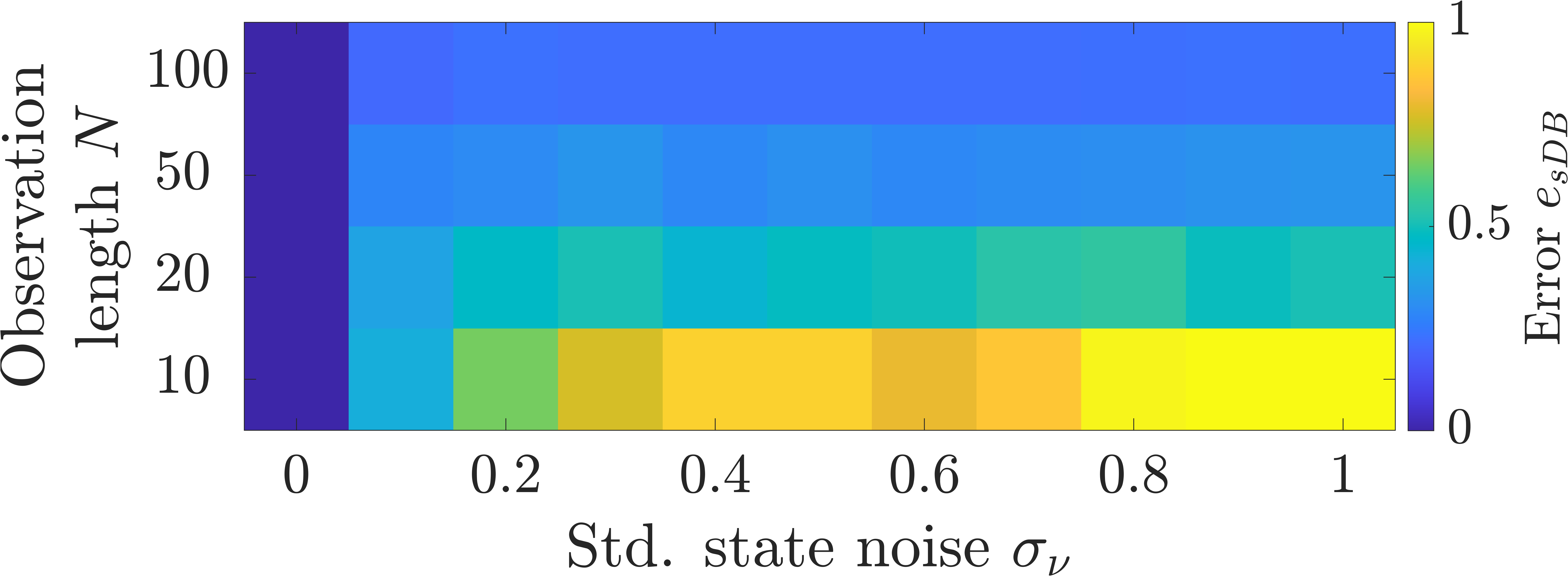}
         \caption{}
         \label{}
     \end{subfigure}
        \caption{Mean relative errors $e = \|\hat{X}-X\|/\|X\|$ for different state noise $\sigma_{\nu}$ and observation length $N$ from simulations of $100$ trajectories for each case. a)~The proposed filter. b)~Simple dead-beat observer. c)~Sliding-window dead-beat observer.}
        \label{fig:Comparison}
\end{figure}

The results are depicted in Fig. \ref{fig:Comparison}. Observe, for $N=10=n$, that all three methods fail to reconstruct the states in a meaningful manner except under zero noise when all three methods provide perfect state estimation (left-most columns in Fig. \ref{fig:Comparison}). Both the $H^*$-filter and dead-beat observer with sliding window are relatively robust against state noise, with the former performing slightly better over almost the whole range of experiments. 
Also, the performance increases with greater $N$, except for the simple dead-beat observer that is unable to recover the initial condition correctly for larger $N$. These results underline our claim that solving for a more complex object via the `overparameterization' approach is in fact beneficial.

\subsection{Example 3: Teaser for the nonlinear case}

To demonstrate a nonlinear case, we revisit the H\'enon map \cite{Mareels2002}. Given is a finite set of observations $y_t,~t=1\cdots N$ from the (assumed unknown) model
\begin{equation}
\begin{array}{lcll}
   x_{1,t+1} &=& {\Theta}_1 + x_{2,t} + {\Theta}_2 x_{1,t}^2, \\
   x_{2,t+1} &=& {\Theta}_{3} x_{1,t}, \\
   y_t &=& (1+\mu_t)x_{1,t},
   \end{array}  \label{henon}
\end{equation}
with states $x_{1,t}, x_{2,t}$ and parameters $\Theta = [\Theta_1,\Theta_2, \Theta_{3}]^T$. The measurements are contaminated with multiplicative, i.i.d. uniform noise $\mu_t \sim U[-\sigma_{\mu},\sigma_{\mu}]$.

Using a little prior knowledge about (\ref{henon}), we propose a class of autoregressive, second order polynomial dynamics (this class is a natural, low complexity case that contains the H\'enon map) with $x_t = [y_t, y_{t-1}]^T$ as state variable:
\begin{equation}
\begin{array}{lcll}
\label{henontakensclass}
    y_{t+1} = \theta_{1,t} + \theta_{2,t}y_t+\theta_{3,t}y_{t-1}
     + \theta_{4,t}y_t y_{t-1}    + \theta_{5,t}y^2_t + \theta_{6,t}y^2_{t-1}.
\end{array}
\end{equation}
We want to simultaneously reconstruct the state $x_t$ (like de-noising in this case) and the parameters, i.e., $\theta_{1,t} = 1, \theta_{2,t} = 0, \theta_{4,t} = 0, \theta_{6,t} = 0$ and, considering the chaotic attractor, $\theta_{3,t} = 0.3$ and $\theta_{5,t} = -1.4, \forall t$. 

The respective state transition $f$ in (\ref{augmentedstatefunction}) can then be formulated with trivial dynamics for the parameters $\theta_t = (\theta_{1,t},\theta_{2,t},\theta_{3,t},\theta_{4,t},\theta_{5,t},\theta_{6,t})^T$. With the state $Z= (y_{1},\theta_{1,1},\theta_{2,1},\theta_{3,1} \cdots \theta_{4,N},\theta_{5,N},\theta_{6,N})^T$, we optimize
\begin{equation}
\label{henonloss}
    \ell(\hat{Z}) = \sum_{t=2}^{N-1} \| \hat{y}_{t+1}-f(\hat{y}_t,\hat{y}_{t-1},\hat{\theta}_t)\|_2^2 +\rho \sum_{t=1}^{N}\|y_t-\hat{y}_{t}\|^2_2 
    +\sum_{t=1}^{N-1} \|\hat{\theta}_{t+1} -\hat{\theta}_{t}\|^2_2 + \lambda \sum_{t=1}^{N} \|\hat{\theta}_{t}\|_1,
\end{equation}
adding a $l1$-regularization term weighted by a scalar $\lambda$ to promote sparsity given the dynamical candidate library in (\ref{henontakensclass}) (see \cite{Rudy2018}, presenting a similar approach) and serving as further measure against `overfitting' noise.

\begin{figure}[]
     \centering
     \includegraphics[width=0.55\textwidth]{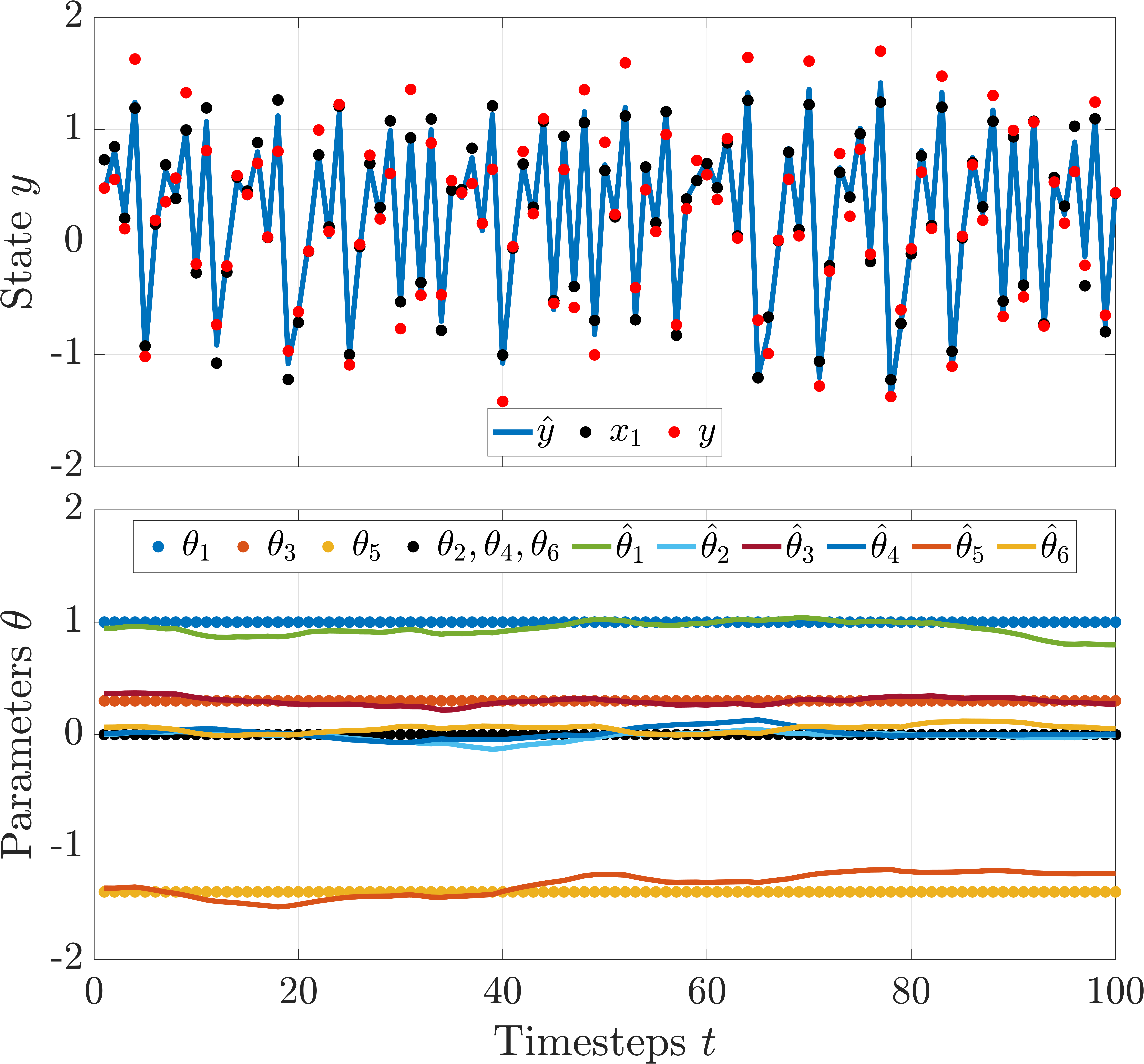}
     \caption{H\'enon map state, output and parameter reconstruction based on the autoregressive polynomial model class (\ref{henontakensclass}). System simulated with $N=100$, random initial condition and uniform noise $\mu_t \sim U[-0.5,0.5]$. Reconstruction with $\rho = 0.1$ and $\lambda = 0.001$.}
     \label{henonexample}
\end{figure}

To minimize (\ref{henonloss}), we employ a single layer neural network as in (\ref{network}) with gradient descent, random NTK initialization and a learning rate $\eta = 0.05$ in \textit{Tensorflow} \cite{Abadi2016TensorFlow:Learning}. Now, given the state-space extension, the output of the network is simply the estimate $\hat{Z} \in R^{7N}$. One particular example of state reconstruction for $\sigma_{\mu} = 0.5$ is displayed in Fig. \ref{henonexample}.
The results of $10$ state trajectories for different noise levels, simulated from (\ref{henon}) with the true $\Theta$ values, and random, but non-divergent initial conditions $x_{1,1}, x_{2,1}$, both sampled from the interval $[0,1]$, are shown in Table \ref{example3} for the case $N=100$, $\rho=0.1, \lambda = 0.001$. Here, we again consider the relative error $\|\hat{X}-X\|/\|X\|$ and $\|\hat{\Theta}-\Theta\|/\|\Theta\|$, with $X = (x_{1,1} \cdots x_{1,N})^T$, $ \hat{X} = (\hat{y}_{1} \cdots \hat{y}_{N})^T$, and $\hat{\Theta}$ being the mean of each estimated $\hat{\theta}_t$ trajectory. 
Importantly, the results in Table \ref{example3} are insensitive to the choice of learning rate and initialization, which is in line with our claim that there are in fact unexpected benefits from an highly `overparameterized' loss that can be minimized by artificial neural networks. That is, transferring the optimization into a high-dimensional space provides ways to circumvent problems associated with classical approaches that recover only an initial condition, which is an optimization of a prohibitively complex loss for systems like the H\'enon map (e.g., \cite{Voss2004}). We hypothesize that the presented overparameterization has, in fact, similar effects as highly overparameterized neural networks in finding global minima \cite{Liu2020a} and avoiding local minima.

\begin{table}
\caption{Results of the state and parameter construction of the H\'enon map for $\rho = 0.1$ and $\lambda = 0.001$.}
\label{example3}
\begin{center}
\begin{tabular}{|c||c| c | c | c | c |}
\hline
Noise $\sigma_{\mu}$ & 0 & 0.1 & 0.2 & 0.5 & 1\\
\hline \hline
Relative error in $X$ & 0.014 & 0.026 & 0.050 & 0.137 & 0.486\\
\hline
Relative error in $\Theta$ & 0.025 & 0.029 & 0.037 & 0.074 & 0.489\\
\hline
\end{tabular}
\end{center}
\end{table}

\section{Conclusions}
 We revisited the state (parameter) estimation problem for a given model class when finitely many noisy measurements are available. In typical batch mode, we propose to estimate the state trajectory as the minimizer of a quadratic loss function defined by the available prior information. In the case of linear dynamics, this estimator has an analytical expression. At each time instant, the state estimate is a linear function of the available measurements. The `overparameterization' provides robustness against noise, and can be tuned to reflect the relative size of state and measurement noise (which can be estimated through the same procedure). Our approach can be extended to address nonlinear problems, such as the joint parameter and state estimation of potentially complex dynamics. The loss function in this case has to be optimized using numerical methods. We argued that a neural network approach to accomplish this optimization provides a way forward. An example to illustrate this has been presented. The analysis is work in progress.

\bibliographystyle{unsrt}  
\bibliography{references}

\end{document}